\documentclass[preprint,sort&compress,final,11pt]{elsarticle}

\usepackage{amsmath}
\usepackage{amssymb}
\usepackage{graphicx}
\usepackage{latexsym}
\usepackage{epstopdf}
\usepackage{natbib}
\usepackage{color}
\usepackage{hyperref}

\newtheorem{theorem}{Theorem}[section]

\newtheorem{lemma}[theorem]{Lemma}

\newtheorem{remark}[theorem]{Remark}
\numberwithin{equation}{section}

\def\Proof{\noindent{\bf Proof.}~}
\def\qed{\hfill$\square$\smallskip}

\journal{\empty}
\date{}

\begin{document}

\begin{frontmatter}

\title{Periodic solutions of semilinear Duffing equations with impulsive effects}

\author[au1]{Yanmin Niu}

\ead[au1]{nym1216@163.com}

\author[au1]{Xiong Li\footnote{ Partially supported by the NSFC (11571041)and the Fundamental Research Funds for the Central Universities. Corresponding author.}}

\address[au1]{School of Mathematical Sciences, Beijing Normal University, Beijing 100875, P.R. China.}

\ead[au1]{xli@bnu.edu.cn}

\begin{abstract}
In this paper we are concerned with the existence of periodic solutions for semilinear Duffing equations with impulsive effects. Firstly for the autonomous one, basing on Poincar\'{e}-Birkhoff twist theorem, we prove the existence of infinitely many periodic solutions. Secondly, as for the nonautonomous case, the impulse brings us great challenges for the study, and there are only finitely many periodic solutions, which is quite different from the corresponding equation without impulses. Here, taking the autonomous one as an auxiliary equation, we find the relation between these two equations and then obtain the result also by Poincar\'{e}-Birkhoff twist theorem.
\end{abstract}

\begin{keyword}
Impulsive differential equations, Periodic solutions, Poincar\'{e}-Birkhoff twist theorem
\end{keyword}

\end{frontmatter}

\section{Introduction}
We are concerned in this paper with the existence of periodic solutions for the second order impulsive differential equation
\begin{equation}\label{maineq}
\left\{\begin{array}{ll}
 x''+g(x)=p(t),\quad t\neq t_{j};\\[0.2cm]
 \Delta x|_{t=t_{j}}=0,\\[0.2cm]
 \Delta x'|_{t=t_{j}}=-2x'\left(t_{j}-\right),\quad j=\pm1,\ \pm2,\cdots,\\
 \end{array}
 \right.
 \end{equation}
where $0\leq t_{1}<2\pi$, $\Delta x|_{t=t_{j}}=x\left(t_{j}+\right)-x\left(t_{j}-\right)$,
$\Delta x'|_{t=t_{j}}=x'\left(t_{j}+\right)-x'\left(t_{j}-\right)$, $g(x),\ p(t)\in C(\mathbb{R},\mathbb{R})$ and $p(t)$ is $2\pi$-periodic. In addition, we assume that the impulsive time is $2\pi$-periodic, that is,
$t_{j+1}=t_{j}+2\pi$ for $j=\pm1,\ \pm2,\cdots$.

This problem comes from Duffing equation
\begin{equation}\label{Duffing}
x''+g(x)=p(t),
\end{equation}
and there is a wide literature dealing with the existence of periodic solutions for the above Duffing equation, not only because of its physical significance, but also the application of various mathematical techniques on it, such as Poincar\'{e}-Birkhoff twist theorem in \cite{Dingsuperquadratic}, \cite{Hartman77}, \cite{Jacobowitz76}, the variational method in \cite{Bahri841}, \cite{Long89} and topological degree or index theories in \cite{Capietto92}, \cite{Capietto90}. Under different assumptions on the function $g$, for example being superlinear, sublinear, semilinear and so on, there are many interesting results on the existence and multiplicity of periodic solutions of \eqref{Duffing}, see \cite{Dingsublinear}, \cite{Pino}, \cite{Qian01}, \cite{Qian052}, \cite{Rebelo} and the references therein. Among these, the existence problem of periodic solutions for semilinear Duffing equations challenges more attentions for its special resonance phenomenon. At resonance, equation \eqref{Duffing} may have no bounded solutions, therefore the crucial point of solving this problem is to exclude the resonance, and there are many studies on it, see \cite{A. C. Lazer69}, \cite{ding1982}, \cite{TRDing82} and \cite{Reissing75}.

Recently, as impulsive equations widely arise in applied mathematics, they attract a lot of attentions and many authors study the general properties of them in \cite{Bainov93}, \cite{Lak}, along with the existence of periodic solutions of impulsive differential equations via fixed point theory in \cite{Nieto97}, \cite{Nieto02}, topological degree theory in \cite{Dong01}, \cite{Qian05}, and the variational method in \cite{Nieto09}, \cite{Sun11}. However, different from the extensive study for second order differential equations without impulsive terms, there are only a few results on the existence and multiplicity of periodic solutions for second order impulsive differential equations.

In \cite{Qian15}, Qian etc considered the superlinear impulsive differential equation
\begin{equation}\label{qianeq}
\left\{
\begin{array}{ll}
 x''+g(x)=p(t,x,x'),\quad \quad\quad\ \  t\neq t_{j};\\[0.2cm]
\Delta x|_{t=t_{j}}=I_{j} (x\left(t_{j}-\right),x'\left(t_{j}-\right)),\\[0.2cm]
\Delta x'|_{t=t_{j}}=J_{j} (x\left(t_{j}-\right),x'\left(t_{j}-\right)),\ \   j=\pm1,\ \pm2,\cdots,
 \end{array}
 \right.
\end{equation}
where $0\leq t_{1}<\cdots<t_{k}<2\pi$, $I_{j}, J_{j}:\mathbb{R}\times\mathbb{R}\rightarrow\mathbb{R}$ are continuous maps, $t_{j+k}=t_{j}+2\pi$ for $j=\pm1,\ \pm2,\cdots$, and $g$ is a continuous function with the superlinear growth condition
$$(g_{0}'):\ \displaystyle\lim_{\mid x\mid\to+\infty}\frac{g(x)}{x}=+\infty.$$ The authors proved via Poincar\'{e}-Birkhoff twist theorem the existence of infinitely many periodic solutions of \eqref{qianeq} with $p=p(t)$, and also the existence
of periodic solutions for non-conservative case with degenerate impulsive terms by developing a new twist fixed point theorem.

In \cite{niu}, we discussed the existence of periodic solutions for the sublinear impulsive differential equation
\begin{equation}\label{niueq}
\left\{\begin{array}{ll}
 x''+g(x)=p(t,x,x'),\quad t\neq t_{j};\\[0.2cm]
 \Delta x|_{t=t_{j}}=a x\left(t_{j}-\right),\\[0.2cm]
 \Delta x'|_{t=t_{j}}=a x'\left(t_{j}-\right),\quad\ \ \ \ \ j=\pm1,\ \pm2,\cdots,\\
 \end{array}
 \right.
 \end{equation}
where $0\leq t_{1}<\cdots<t_{k}<2\pi$, $a>0$ is a constant, $t_{j+k}=t_{j}+2\pi$ for $j=\pm1,\ \pm2,\cdots$, and $g$ is a continuous function with the sublinear growth condition
$$(g_{0}):\ \displaystyle\lim_{\mid x\mid\to+\infty}\frac{g(x)}{x}=0.$$
Basing on the Poincar\'{e}-Bohl fixed point theorem and a new twist fixed point theorem recently established by Qian etc in \cite{Qian15}, we obtained the existence of harmonic solutions and subharmonic solutions, respectively. The impulsive function is especially chosen to keeping the arguments of trajectories unchanged under the polar coordinate, such that the obstacle causing by the impulse can be solved.

In this article, we discuss semilinear Duffing equation \eqref{maineq} with the impulse, which is different from the superlinear or sublinear case and there is few papers studying on it up to now. As we all know, the existence of impulses, even the simplest impulsive function, may cause complicated dynamic phenomena and bring great difficulties to study. The behavior of solutions with impulsive effects may have great changes compared with solutions without impulses. Here we take the simplest linear equation as an example and consider
\begin{equation}\label{example equation}
x''+x=0,
\end{equation}
with the impulsive condition
\begin{equation}\label{example impulse}
 \Delta x|_{t=t_{j}}=0,\quad\quad\ \Delta x'|_{t=t_{j}}=-2x'\left(t_{j}-\right),\quad j=\pm1,\ \pm2,\cdots,
\end{equation}
where $0<t_{1}<2\pi$, $t_{j+1}=t_{j}+2\pi$ for $j=\pm1,\ \pm2,\cdots$.
It is easy to see that without the impulse, all solutions of \eqref{example equation} are $2\pi$-periodic and  has  the form
$$
x(t)=C_{1}\cos t+C_{2}\sin t,
$$
where $C_1$ and $C_2$ are arbitrary constants. However, with $t_{1}=\frac{\pi}{2}$ in \eqref{example impulse}, the solution starting from $(x(0),x'(0))=(x_{0},y_{0})$ has the form
\begin{equation}\label{solution form}
x(t)=\left\{\begin{array}{lll}
x_{0}\cos t+y_{0}\sin t,\quad\quad & t\in[0,\frac{\pi}{2}),\\[0.2cm]
y_{0}\cos (t-\frac{\pi}{2})+x_{0}\sin (t-\frac{\pi}{2}),\quad\quad &t\in[\frac{\pi}{2},\frac{5\pi}{2}),\\[0.2cm]
y_{0}\cos (t-\frac{5\pi}{2})-x_{0}\sin (t-\frac{5\pi}{2}),\quad\quad &t\in[\frac{5\pi}{2},4\pi],\\
\end{array}
\right.
\end{equation}
from which one has $(x(2\pi),x'(2\pi))=(-x_{0},y_{0})$ and $(x(4\pi),x'(4\pi))=(x_{0},y_{0})$.
Therefore for $x_{0}\neq 0$, the solutions of \eqref{example equation} have the least period $4\pi$, which implies that the impulse breaks the rule of periods of solutions. On the other hand, if $x_{0}=0$,
$$
x(t)=y_{0}\sin t,\quad\quad t\in \mathbb{R},
$$
the period of which is $2\pi$. In this case, the impulse exists in name only and has no any effects on solutions. More specific discussions about this linear impulsive equation are in \cite{Bainov93} Chapter II.

The construction of the impulsive function in \eqref{maineq} is also inspired by the oscillator with obstacle, which has the form
\begin{equation}\label{oscillator}
\left\{\begin{array}{ll}
 x''+g(x)=f(t),\ \ x(t)> 0;\\[0.2cm]
 x(t)\geq 0;\\[0.2cm]
 x(t_{0})=0\ \Rightarrow x'(t_{0}+)=-x'(t_{0}-),\\
 \end{array}
 \right.
 \end{equation}
where $g$ and $f$ are continuous functions.
This equation has actual physical backgrounds in classical billiard, in electric and magnetic fields, three body problem and so on, see \cite{pboyland}, \cite{corbera}, \cite{hlamba}. Equation \eqref{oscillator} can be thought of the model of the motion of a particle which is attached to a spring that pushes the particle against a rigid wall situated at $x=0$, and at this barrier the particle bounces elastically. At the bouncing time $t_{0}$ satisfying $x(t_{0})=0$, the velocity change $x'(t_{0}+)=-x'(t_{0}-)$ is same as the impulsive function in \eqref{maineq}. In both cases, the motions  are actually continuous and just the velocity changes the sign at the impulsive/bouncing times. However, in \eqref{maineq} the impulsive times are given in advance while in \eqref{oscillator} the bouncing times depend on the specific solutions, even any positive solutions never meet the wall.

It is well known that the solutions of general autonomous Duffing equations are a family of closed curves in the phase plane. Due to the existence of impulses, the closed curves may not be preserved and even be broken. For this point, the impulsive function in \eqref{maineq} guarantees the motion of trajectories still on the closed curves, which reduces the difficulty causing by impulses.

The rest is organized as follows. The main results about the existence of periodic solutions of \eqref{maineq} and \eqref{maineqautonomic}, together with some basic lemmas and symbols, are given in Section 2. In Section 3, we consider the autonomous equation with the same impulsive condition as \eqref{maineq}, then obtain infinitely many $2\pi$-periodic solutions by Poincar\'{e}-Birkhoff twist theorem. In the last section, from the discussion about the autonomous equation, the relation between arguments of the two equations under the polar coordinate is established. Again it follows from Poincar\'{e}-Birkhoff twist theorem that there are finitely many $2\pi$-periodic solutions of \eqref{maineq}.

\section{The main results}

For convenience, we introduce some basic lemmas and symbols which are used in the sequel.
Firstly consider Duffing equation without impulses and perturbation
\begin{equation}\label{Duffing eq basic}
x''+g(x)=0,
\end{equation}
and its equivalent system
\begin{equation}\label{Duffing eq basic polar}
x'=y,\quad y'=-g(x).
\end{equation}
This is a planar autonomous system whose orbits are curves determined by the following equation
\begin{equation}\label{Vfunction}
V(x,y)=\frac{1}{2}y^{2}+G(x)=c,
\end{equation}
where $G(x)=\int_{0}^{x}g(u)du$ and $c>0$ is a parameter. Then $V^{-1}(c)$ is star-shaped, the proof of which we omit and one can refer to \cite{ding1982}.

\begin{lemma}\label{starshape}(\cite{ding1982})
If the condition $(H_{2})$ below holds, then $V^{-1}(c)$ is a closed curve for $c\geq c_{0}$ which is star-shaped about the origin, where $c_0$ is a positive constant.
\end{lemma}

Denote the curve $V^{-1}(c)$ by $\Gamma_{c}$. It follows from Lemma \ref{starshape} that each curve $\Gamma_{c}\ (c\geq c_{0})$ intersects the $x$-axis at two points: $(h(c),0)$ and $(-h_{1}(c),0)$, where $h(c)>0$ and $h_{1}(c)>0$ are uniquely determined by the formula
\begin{equation}\label{definition of h}
G(h(c))=G(-h_{1}(c))=c.
\end{equation}
Let $(x(t),y(t))$ be any solutions of \eqref{Duffing eq basic polar} whose orbit is $\Gamma_{c}\ (c\geq c_{0})$. Obviously, this solution is periodic and denote by $\tau(c)$ the least positive period of it. It can be induced from equation \eqref{Vfunction} that
$$\tau(c)=\sqrt{2}\displaystyle\int_{-h_{1}(c)}^{h(c)}\frac{du}{\sqrt{c-G(u)}}.$$

Generally speaking,  Poincar\'{e}-Birkhoff twist theorem is a powerful tool to obtain the existence of periodic solutions. Here we briefly restate a generalized form of it in \cite{dingwy}.

Let $\mathcal{D}$ denote an annular region in the $(x,y)$-plane. The boundary of $\mathcal{D}$ consists of two simple closed curves: the inner boundary curve $C_{1}$ and the outer boundary curve $C_{2}$. Let $\mathcal{D}_{1}$ denote the simple connected open set bounded by $C_{1}$. Consider an area-preserving mapping $\mathrm{T}: \mathbb{R}^{2}\rightarrow \mathbb{R}^{2}$. Suppose that $\mathrm{T}(\mathcal{D})\subset \mathbb{R}^{2}-{0}$, where $0$ is the origin. Suppose $(\gamma,\theta)$ be the polar coordinate of $(x,y)$, that is, $x=\gamma\cos\theta, y=\gamma\sin\theta$. Assume the restriction $T|\mathcal{D}$ is given by
$$
\gamma^{\ast}=f(\gamma,\theta),\quad \quad \theta^{\ast}=\theta+g(\gamma,\theta),
$$
where $f$ and $g$ are continuous in $(\gamma,\theta)$ and $2\pi$-periodic in $\theta$.
\begin{lemma} (\cite{dingwy}) Besides the above mentioned assumptions, we assume that
\begin{enumerate}\label{PBfixed}
\item $C_{1}$ is star-shaped about the origin;
\item $0\in \mathrm{T}(\mathcal{D}_{1})$;
\item $g(\gamma,\theta)>0(<0),\ (\gamma\cos\theta,\gamma\sin\theta)\in C_{1};\\ g(\gamma,\theta)<0(>0),\ (\gamma\cos\theta,\gamma\sin\theta)\in C_{2}$.
\end{enumerate}
Then $\mathrm{T}$ has at least two fixed points in $\mathcal{D}$.
\end{lemma}

Before stating the main results, we give the following hypotheses:

$(H_{1})$ Let $g(x)\in C^1(\mathbb{R},\mathbb{R})$, and $K$ be a positive constant, such that
$$|g'(x)|\leq K, \quad\quad x\in \mathbb{R};$$

$(H_{2})$ There exist two constants $A_{0}>0$ and $M_{0}>0$, such that
$$x^{-1}g(x)\geq A_{0}, \quad\quad |x|\geq M_{0}.$$

$(H_{3})$ There exist a constant $\alpha>0$, an integer $m>0$, and two sequences ${a_{k}}$ and ${b_{k}}$, such that $a_{k}\rightarrow \infty$ and $b_{k}\rightarrow \infty$ as $k\rightarrow \infty$; and moreover\\
$$\tau(a_{k})<\frac{2\pi}{m}-\alpha,\quad \quad \tau(b_{k})>\frac{2\pi}{m}+\alpha.$$

These assumptions are completely same as that in \cite{ding1982}, and Ding in this paper constructed a function $g$ satisfying $(H_{1})-(H_{3})$ very skillfully.

Our results on equation \eqref{maineq} are different for $p(t)=0$ and $p(t)\neq0$. Then we will discuss the two cases respectively in the following parts. The autonomous impulsive equation has the form
\begin{equation}\label{maineqautonomic}
\left\{\begin{array}{ll}
 z''+g(z)=0,\quad t\neq t_{j};\\[0.2cm]
 \Delta z|_{t=t_{j}}=0,\\[0.2cm]
 \Delta z'|_{t=t_{j}}=-2z'\left(t_{j}-\right),\quad\ \ \ \ \ j=\pm1,\ \pm2,\cdots,\\
 \end{array}
 \right.
 \end{equation}
where $0\leq t_{1}<2\pi$, $t_{j+1}=t_{j}+2\pi$ for $j=\pm1,\ \pm2,\cdots$.

\begin{theorem}\label{mainresultautonomic}
Assume that $(H_{1})-(H_{3})$ hold. Then equation \eqref{maineqautonomic} has infinitely many $2\pi$-periodic solutions.
\end{theorem}

As for \eqref{maineq}, due to the existence of the impulse, we only guarantee that the number of $2\pi$-periodic solutions is finite, which is not totally same as the result on semilinear Duffing equation without the impulse in \cite{ding1982}.

\begin{theorem}\label{mainresult}
Assume that $(H_{1})-(H_{3})$ hold and $p(t)$ is a $2\pi$-periodic continuous function. Then equation \eqref{maineq} has finitely many $2\pi$-periodic solutions.
\end{theorem}

\section{Autonomous Duffing impulsive equation}
 We start with equation \eqref{maineqautonomic}, which is a special case of \eqref{maineq} with $p(t)=0$, and its equivalent system is
\begin{equation}\label{auxiliarytwoeq}
\left\{\begin{array}{ll}
 z'=w,\quad w'=-g(z),\quad t\neq t_{j};\\[0.2cm]
 \Delta z|_{t=t_{j}}=0,\\[0.2cm]
 \Delta w|_{t=t_{j}}=-2w\left(t_{j}-\right),\quad\ \ \ \ \ j=\pm1,\ \pm2,\cdots.\\
 \end{array}
 \right.
 \end{equation}
Let $\left(\overline{z}(t,z,w),\overline{w}(t,z,w)\right)$ be the solution of \eqref{auxiliarytwoeq} with the initial point $\left(\overline{z}(0),\overline{w}(0)\right)=(z,w)$. It is not hard to show that every such solution exists on the whole $t$-axis under the condition $(H_{1})$ (see \cite{Bainov93}). Then the Poincar\'{e} map $\mathrm{P_{1}}: \mathbb{R}^{2}\rightarrow \mathbb{R}^{2}$ is well defined by
$$
(z,w)\rightarrow\left(\overline{z}(2\pi,z,w),\overline{w}(2\pi,z,w)\right).
$$

Next, we take the transform $\overline{z}(t)=\rho(t)\cos\varphi(t),\ \overline{w}(t)=\rho(t)\sin\varphi(t)$ for system \eqref{auxiliarytwoeq}. Then the resulting equations for $\rho(t)$ and $\varphi(t)$ are
\begin{equation}\label{auxiliarytwopolareq}
\left\{\begin{array}{ll}
 \rho'=\rho\cos\varphi\cdot\sin\varphi-g(\rho\cos\varphi)\sin\varphi,\\[0.2cm]
 \varphi'=-\sin^{2}\varphi-\displaystyle\frac{1}{\rho}g(\rho\cos\varphi)\cos\varphi,\quad t\neq t_{j};\\[0.4cm]
 \rho(t_{j}+)=\rho(t_{j}-),\\[0.2cm]
 \varphi(t_{j}+)=\arctan(-\tan\varphi (t_{j}-)),\ j=\pm1,\ \pm2,\cdots.\\
 \end{array}
 \right.
 \end{equation}
Let $\left(\overline{\rho}(t,\rho,\varphi),\overline{\varphi}(t,\rho,\varphi)\right)$ be the solution of \eqref{auxiliarytwopolareq} with the initial point $\left(\overline{\rho}(0),\overline{\varphi}(0)\right)=(\rho,\varphi)$. Now we estimate the variation of $\overline{\varphi}(t,\rho,\varphi)$ during the time interval $[0,2\pi]$.

\begin{lemma}\label{twistlemma}
Let $\Phi(\rho,\varphi)=\overline{\varphi}(2\pi,\rho,\varphi)-\varphi$. Then for $k$ sufficiently large, there exist two positive constants $\beta_{1},\ \beta_{2}$, such that
\begin{equation}\label{twistequaiton}
\left\{
\begin{array}{ll}
\Phi(\rho,\varphi)\leq-2m\pi-\beta_{1},\quad\quad(\rho\cos\varphi,\rho\sin\varphi)\in \Gamma_{a_{k}};\\[0.2cm]
\Phi(\rho,\varphi)\geq-2m\pi+\beta_{2},\quad\quad(\rho\cos\varphi,\rho\sin\varphi)\in \Gamma_{b_{k}},\\
\end{array}
 \right.
 \end{equation}
where $m,\ a_{k},\ b_{k}$ are given in $(H_{3})$.
\end{lemma}

\Proof Suppose $(z,w)=(\rho\cos\varphi,\rho\sin\varphi)\in \Gamma_{a_{k}}$. Consider the solution $\left(\overline{\rho}(t,\rho,\varphi),\overline{\varphi}(t,\rho,\varphi)\right)$ of \eqref{auxiliarytwopolareq}. Without loss of generality, let the initial point $(\rho,\varphi)$ be on the positive $x$-axis, that is, $(\rho,\varphi)=(\rho,0)$. From $(H_{2})$, there exist two constants $c_{0}>0$ and $A_{1}>0$ such that if $c\geq c_{0}$, then
\begin{equation}\label{A1}
w^{2}+zg(z)\geq A_{1}(z^{2}+w^{2}), \quad \quad (z,w)\in \Gamma_{c}.
\end{equation}
It follows from \eqref{A1} and the second equation of \eqref{auxiliarytwopolareq} that for $k$ large enough satisfying $a_{k}\geq c_{0}$,
\begin{equation}\label{derivativeofvarphi}
 \overline{\varphi}'(t,\rho,\varphi)\leq -A_{1}.
\end{equation}

Since the solution of \eqref{Duffing eq basic polar} lying in $\Gamma_{a_k}$ has the least period $\tau(a_{k})$, we know that the time in which $\overline{\varphi}(t)$ has a decrement $2\pi$ without the impulse is just $\tau(a_{k})$. Denote
$$
\Phi(\rho,\varphi)=\overline{\varphi}(2\pi,\rho,\varphi)-\varphi=\overline{\varphi}(2\pi)-\overline{\varphi}(0)=-2l\pi-\sigma,
$$
where $l\geq 0$ is an integer, and $0\leq\sigma<2\pi$. Let $t_{\sigma}$ denote the time in which $\overline{\varphi}(t)$ decreases from $\varphi-2l\pi$ to $\varphi-2l\pi-\sigma$. Then there are three possible cases according to different positions of $\overline{\varphi}(t)$ when the impulse occurs.

{\bf Case i}\, The impulsive time $t_{1}$ occurs during the time interval in which $\overline{\varphi}(t)$ travels the front $l-1$ loops, that is, in the time interval $[t_{1},2\pi]$ , $\overline{\varphi}(t)$ can travel at least one loop. Then there exists $0<\tau^{\ast}\leq\tau(a_{k})$ such that
$$(l-1)\tau(a_{k})+\tau^{\ast}+t_{\sigma}=2\pi,$$
and the impulsive time $t_{1}$ is out of the time interval in which $\overline{\varphi}(t)$ decreases from $\varphi-2l\pi$ to $\varphi-2l\pi-\sigma$. Assume the plane coordinate of the solution at $t_{1}\pm$ be $(z(t_{1}\pm),w(t_{1}\pm))$, then $\tau^{\ast}$ stands for the time of $\overline{\varphi}(t)$ traveling the argument of $2\pi$ under the impulsive effect and has the following form
$$
\tau^{\ast}=\tau(a_{k})-\sqrt{2}\displaystyle\int_{z(t_{1}-)}^{h(a_{k})}\frac{du}{\sqrt{a_{k}-G(u)}},
$$
when $(z(t_{1}-),w(t_{1}-))$ is in the first or the second quadrant, and has the form
$$
\tau^{\ast}=\tau(a_{k})-\sqrt{2}\displaystyle\int_{-h_{1}(a_{k})}^{z(t_{1}+)}\frac{du}{\sqrt{a_{k}-G(u)}},
$$
when $(z(t_{1}-),w(t_{1}-))$ is in the third or the forth quadrant, where $h,\ h_{1}$ are defined in \eqref{definition of h}.

In this case we infer $l\geq1$. Since $0\leq t_{\sigma}<\tau(a_{k})$, by the condition $(H_{3})$ it holds that
$$
\begin{array}{lll}
2\pi&=&(l-1)\tau(a_{k})+\tau^{\ast}+t_{\sigma}\\[0.2cm]
&<&(l-1)\tau(a_{k})+\tau(a_{k})+\tau(a_{k})\\[0.2cm]
&=&(l+1)\tau(a_{k})\\[0.2cm]
&\leq&(l+1)\left(\frac{2\pi}{m}-\alpha\right),
\end{array}
$$
which implies $l\geq m$. If $l\geq m+1$, we have
\begin{equation}\label{m+1}
\Phi(\rho,\varphi)=-2l\pi-\sigma\leq-2l\pi\leq-2(m+1)\pi.
\end{equation}
Now, assume $l=m$. Then we have
\begin{equation}\label{estimate the time not equal to oneloop}
\begin{array}{lll}
t_{\sigma}&=&2\pi-(m-1)\tau(a_{k})-\tau^{\ast}\\[0.2cm]
&>&2\pi-(m-1)\tau(a_{k})-\tau(a_{k})\\[0.2cm]
&=&2\pi-m\tau(a_{k})\\[0.2cm]
&\geq&2\pi-m\left(\frac{2\pi}{m}-\alpha\right)=m\alpha.
\end{array}
\end{equation}
By \eqref{derivativeofvarphi} and \eqref{estimate the time not equal to oneloop}, we obtain
$$
-\sigma=\displaystyle\int_{(l-1)\tau(a_{k})+\tau^{\ast}}^{(l-1)\tau(a_{k})+\tau^{\ast}+t_{\sigma}} \overline{\varphi}'(t,\rho,\varphi)dt\leq-A_{1}t_{\sigma}\leq-m\alpha A_{1}.
$$
Thus
\begin{equation}\label{m}
\Phi(\rho,\varphi)=-2l\pi-\sigma\leq-2m\pi-m\alpha A_{1}.
\end{equation}
Combining \eqref{m+1} and \eqref{m} yields the first inequality of \eqref{twistequaiton}, where $\beta_{1}=\min\{2\pi,m\alpha A_{1}\}$.

{\bf Case ii}\, The impulsive time $t_{1}$ occurs at the $l$th loop and in the time interval $[t_{1},2\pi]$, $\overline{\varphi}(t)$ travels no more that one loop. Then there exists $0<\tau^{\ast}\leq\tau(a_{k})$ such that
$$
(l-1)\tau(a_{k})+\tau^{\ast}+t_{\sigma}=2\pi,
$$
and the impulsive time $t_{1}$ has impact on the time interval in which $\overline{\varphi}(t)$ decreases from $\varphi-2l\pi$ to $\varphi-2l\pi-\sigma$.
More specifically, this case happens when $(z(t_{1}-),w(t_{1}-))$ is in the first quadrant or the second quadrant, such that $\tau^{\ast}$ has the form
$$
\tau^{\ast}=\tau(a_{k})-\frac{\sqrt{2}}{2}\displaystyle\int_{z(t_{1}-)}^{h(a_{k})}\frac{du}{\sqrt{a_{k}-G(u)}}.
$$

Also we infer $l\geq1$ and similarly to the proof of case $(i)$, we obtain $l\geq m$, together with \eqref{m+1} and \eqref{estimate the time not equal to oneloop}. By \eqref{derivativeofvarphi} and \eqref{estimate the time not equal to oneloop}, we have
$$
\begin{array}{lll}
-\sigma&=&\frac{1}{2}\left(\overline{\varphi}'(t_{1}+,\rho,\varphi)-\overline{\varphi}'(t_{1}+,\rho,\varphi)\right)+\displaystyle\int_{(l-1)\tau(a_{k})+\tau^{\ast}}^{(l-1)\tau(a_{k})+\tau^{\ast}+t_{\sigma}} \overline{\varphi}'(t,\rho,\varphi)dt\\[0.2cm]
&\leq&\displaystyle\int_{(l-1)\tau(a_{k})+\tau^{\ast}}^{(l-1)\tau(a_{k})+\tau^{\ast}+t_{\sigma}}\overline{\varphi}'(t,\rho,\varphi)dt\\[0.4cm]
&\leq&-A_{1}t_{\sigma}\\[0.2cm]
&\leq&-m\alpha A_{1}.
\end{array}
$$
Thus \eqref{m} holds. Again combining \eqref{m+1} and \eqref{m} yields the first inequality of \eqref{twistequaiton}, where $\beta_{1}=\min\{2\pi,m\alpha A_{1}\}$.

{\bf Case iii}\, The movement of $\overline{\varphi}(t)$ during the front $l$ loops is out of the effect of the impulsive time $t_{1}$. Then $$l\tau(a_{k})+t_{\sigma}=2\pi$$ and the impulsive time $t_{1}$ is in the time interval in which $\overline{\varphi}(t)$ decreases from $\varphi-2l\pi$ to $\varphi-2l\pi-\sigma$. In this case,
$$
t_{1}\in\left(l\cdot\tau(a_{k}),2\pi\right).
$$
Since $0\leq t_{\sigma}<\tau(a_{k})$, by the condition $(H_{3})$ it holds that
$$
2\pi=l\cdot\tau(a_{k})+t_{\sigma}<(l+1)\tau(a_{k})\leq(l+1)\left(\frac{2\pi}{m}-\alpha\right),
$$
which implies $l\geq m$. If $l\geq m+1$, \eqref{m+1} holds.
Assume $l=m$. Then we have
\begin{equation}\label{estimate the time not equal to oneloop case3}
t_{\sigma}=2\pi-m\tau(a_{k})\geq2\pi-m\left(\frac{2\pi}{m}-\alpha\right)=m\alpha.
\end{equation}
By \eqref{derivativeofvarphi} and \eqref{estimate the time not equal to oneloop case3}, we obtain
$$
\begin{array}{lll}
-\sigma&=&\displaystyle\int_{l\cdot\tau(a_{k})}^{t_{1}-} \overline{\varphi}'(t,\rho,\varphi)dt+
\left(\overline{\varphi}'(t_{1}+,\rho,\varphi)-\overline{\varphi}'(t_{1}-,\rho,\varphi)\right)+\displaystyle\int_{t_{1}+}^{2\pi} \overline{\varphi}'(t,\rho,\varphi)dt\\[0.4cm]
&\leq&-A_{1}(t_{1}-l\cdot\tau(a_{k}))-A_{1}(2\pi-t_{1})\\[0.2cm]
&=&-A_{1}t_{\sigma}\\[0.2cm]
&\leq&-m\alpha A_{1}.
\end{array}
$$
Thus
$$
\Phi(\rho,\varphi)=-2l\pi-\sigma\leq-2m\pi-m\alpha A_{1}.
$$
Similarly, we obtain the first inequality of \eqref{twistequaiton}, where $\beta_{1}=\min\{2\pi,m\alpha A_{1}\}$.
%%%%%%%%%%%%%%%%%%%%%%%%%%%%%%%%%%%%%%%%%%%%%%%%%%%%%%%%%%%%%%%%%%%%%%%%%%%%%

We are now in a position to prove the second inequality of \eqref{twistequaiton}. Suppose $(z,w)=(\rho\cos\varphi,\rho\sin\varphi)\in \Gamma_{b_{k}}$, where $k$ is sufficiently large such that $b_{k}\geq c_{0}$. We see that the time in which $\overline{\varphi}(t)$ has a decrement $2\pi$ without the impulse is just $\tau(b_{k})$. Denote
$$
\Phi(\rho,\varphi)=\overline{\varphi}(2\pi,\rho,\varphi)-\varphi=\overline{\varphi}(2\pi)-\overline{\varphi}(0)=-2q\pi+\xi,
$$
where $q\geq 1$ is an integer, and $0<\xi\leq2\pi$. Let $t_{\xi}$ denote the time in which $\overline{\varphi}(t)$ decreases from $\varphi-2q\pi+\xi$ to $\varphi-2q\pi$. In the following proof, we always assume that during the time interval from $2\pi$ to $2\pi+t_{\xi}$, there is under no influence of another impulsive time. In fact, if there is another impulsive time $t_{2}$ in $(2\pi,2\pi+t_{\xi})$, the proof below is similar and we just need to change some signs in \eqref{xi}. Under these assumptions, the three cases discussed above on $\Gamma_{a_{k}}$ can be simplified into one case: among the $q$ loops, only one loop is effected by the impulse. Then there exists $\tau^{\wedge}\in(0,\tau(b_{k})]$ such that
$$
(q-1)\tau(b_{k})+\tau^{\wedge}-t_{\xi}=2\pi,
$$
where
$$
\tau^{\wedge}=\tau(b_{k})-\sqrt{2}\displaystyle\int_{z(t_{1}-)}^{h(a_{k})}\frac{du}{\sqrt{b_{k}-G(u)}},
$$
when $(z(t_{1}-),w(t_{1}-))$ is in the first or the second quadrant, and
$$
\tau^{\wedge}=\tau(b_{k})-\sqrt{2}\displaystyle\int_{-h_{1}(b_{k})}^{z(t_{1}+)}\frac{du}{\sqrt{b_{k}-G(u)}},
$$
when $(z(t_{1}-),w(t_{1}-))$ is in the third or the forth quadrant.

Since $0<\tau^{\wedge}\leq\tau(b_{k})$, let $\tau^{\wedge}=\mu\tau(b_{k})$, then $0<\mu\leq1$.
By the condition $(H_{3})$ and $0<t_{\xi}\leq\tau(b_{k})$, it holds that
\begin{equation}\label{estimate q}
2\pi=(q-1)\tau(b_{k})+\tau^{\wedge}-t_{\xi}\geq(q-2+\mu)\tau(b_{k})>(q-2+\mu)\left(\frac{2\pi}{m}+\alpha\right),
\end{equation}
provided $q>1$. Notice that if $q=1$,
$$
\Phi(\rho,\varphi)=-2\pi+\xi\geq-2m\pi+\xi,
$$
then letting $\beta_{2}=\xi$, we obtain the result. We just consider the case $q>1$. It follows from \eqref{estimate q} that $q\leq m$. Notice that if $m$ given in $(H_{3})$ equals to $1$, then $q=1$. Similarly to the discussion above, we obtain the result.
 Next assume $m>1$. If $q\leq m-1$, we have
\begin{equation}\label{m+1b}
\Phi(\rho,\varphi)=-2q\pi+\xi>-2q\pi=-2(m-1)\pi.
\end{equation}
Now if $q=m$, writing $\tau^{\wedge}=\tau(b_{k})-\eta$, we have
\begin{equation}
\begin{array}{lll}
t_{\xi}&=&-2\pi+(m-1)\tau(b_{k})+\tau^{\wedge}\\[0.2cm]
&=&-2\pi+m\tau(b_{k})-\eta\\[0.2cm]
&\geq&-2\pi+m\left(\frac{2\pi}{m}+\alpha\right)-\eta\\[0.2cm]
&=&m\alpha-\eta.
\end{array}
\end{equation}
Since $\xi>0$, $t_{\xi}>\frac{1}{2}t_{\xi}>0$, then there exists $\varepsilon_{\xi}>0$, such that
\begin{equation}\label{estimate the time not equal to oneloopb}
t_{\xi}\geq\max\{\varepsilon_{\xi},m\alpha-\eta\}\triangleq M_{\xi}>0.
\end{equation}
By \eqref{derivativeofvarphi} and \eqref{estimate the time not equal to oneloopb}, we have
\begin{equation}\label{xi}
\begin{array}{lll}
\xi&=&\displaystyle\int_{(q-1)\tau(b_{k})+\tau^{\wedge}}^{(q-1)\tau(b_{k})+\tau^{\wedge}-t_{\xi}} \overline{\varphi}'(t,\rho,\varphi)dt\\[0.4cm]
&=&-\displaystyle\int_{(q-1)\tau(b_{k})+\tau^{\wedge}-t_{\xi}}^{(q-1)\tau(b_{k})+\tau^{\wedge}} \overline{\varphi}'(t,\rho,\varphi)dt\\[0.4cm]
&\geq&-(-A_{1}t_{\xi})\\[0.2cm]
&\geq&M_{\xi}A_{1}.
\end{array}
\end{equation}
Thus
\begin{equation}\label{mb}
\Phi(\rho,\varphi)=-2q\pi+\xi\geq-2m\pi+M_{\xi}A_{1}.
\end{equation}
Combining \eqref{m+1b} and \eqref{mb} yields the second inequality of \eqref{twistequaiton}, where $\beta_{2}=\min\{2\pi,M_{\xi}A_{1}\}$.

The proof of Lemma \ref{twistlemma} is then completed.\qed

\begin{remark}\label{remark} $\beta_{2}$ is related to $k$. Indeed, for each $k$ sufficiently large, $\beta_{2}=\min\{2\pi,M_{\xi}A_{1}\}>0$, where $M_{\xi}\triangleq\max\{\varepsilon_{\xi},m\alpha-\eta\}>0$ may depend on $k$. However, $\beta_{1}=\min\{2\pi,m\alpha A_{1}\}$ is independent of $k$.
\end{remark}

Under the polar coordinate, the Poincar\'{e} map $\mathrm{P_{1}}$ ca be written in the form
\begin{equation}\label{p mapping form autonomic}
\rho^{\ast}=\overline{\rho}(2\pi,\rho,\varphi),\quad \quad \varphi^{\ast}=\overline{\varphi}(2\pi,\rho,\varphi)+2k\pi,
\end{equation}
where $k$ is an arbitrary integer. It can be easily seen that if $(\rho,\varphi)$ satisfies
\begin{equation}\label{rho positive}
\overline{\rho}(t,\rho,\varphi)>0, \quad \quad t\in[0,2\pi],
\end{equation}
then $\overline{\varphi}(2\pi,\rho,\varphi)$ is well defined and continuous in $(\rho,\varphi)$, moreover,
$$
\overline{\varphi}(2\pi,\rho,\varphi+2\pi)=\overline{\varphi}(2\pi,\rho,\varphi)+2\pi.
$$

Now, let $\Gamma_{a_{k}}$ and $\Gamma_{b_{k}}$ be the curves given by Lemma \ref{starshape}, where the specified parameters $a_{k},\ b_{k}\geq c_{0}$ are given by $(H_{3})$, for $k\geq n_{0}$. We can rearrange ${a_{k}}$ and ${b_{k}}$ such that $a_{k}<b_{k}<a_{k+1}$ for $k\geq n_{0}$. Then $\Gamma_{a_{k}}$ and $\Gamma_{b_{k}}$ bound an annular region $\mathcal{A}_{k}$, and $\Gamma_{b_{k}}$ and $\Gamma_{a_{k+1}}$ bound another annular region $\mathcal{B}_{k}$, for $k\geq n_{0}$.

It is well known that each fixed point of $\mathrm{P_{1}}$ corresponds to a $2\pi$-periodic solution of \eqref{maineqautonomic}. In the sequel, we will apply Lemma \ref{PBfixed} to show that $\mathrm{P_{1}}$ has at least two fixed points in each $\mathcal{A}_{k}$ and $\mathcal{B}_{k}$ for sufficiently large $k$. As a consequence, \eqref{maineqautonomic} has an infinite class of $2\pi$-periodic solutions.

We turn to prove Theorem \ref{mainresultautonomic} by the twist theorem.\\

\noindent {\bf Proof of Theorem \ref{mainresultautonomic}} Assume that $a_{k}, b_{k}\geq c_{0}$ for $k\geq n_{0}$, where $n_{0}$ is large enough and $c_{0}$ is given in Lemma \ref{starshape}. Denote by $\mathcal{A}_{k}$ the region bounded by $\Gamma_{a_{k}}$ and $\Gamma_{b_{k}}$. Thus the restriction $\mathrm{P_{1}}|\mathcal{A}_{k}$ can be written in \eqref{p mapping form autonomic}, where we put the integer $k=m$. Then \eqref{p mapping form autonomic} can be rewritten in the form of
$$
\rho^{\ast}=\overline{\rho}(2\pi,\rho,\varphi),\quad \quad \varphi^{\ast}=\varphi+\Phi_{1}(\rho,\varphi),
$$
with $\Phi_{1}(\rho,\varphi)=\overline{\varphi}(2\pi,\rho,\varphi)-\varphi+2m\pi=\Phi(\rho,\varphi)+2m\pi$. By Lemma \ref{twistlemma}, we obtain
$$
\left\{
\begin{array}{ll}
\Phi_{1}(\rho,\varphi)\leq-\beta_{1}<0,\quad\quad(\rho\cos\varphi,\rho\sin\varphi)\in \Gamma_{a_{k}};\\[0.2cm]
\Phi_{1}(\rho,\varphi)\geq\beta_{2}>0,\quad\quad(\rho\cos\varphi,\rho\sin\varphi)\in \Gamma_{b_{k}},\\
\end{array}
 \right.
$$
for $k\geq n_{0}$.

This proves the validity of condition $3$ of Lemma \ref{PBfixed} for the restriction $\mathrm{P_{1}}|\mathcal{A}_{k}$ ($k\geq n_{0}$). Since $k$ can be chosen large enough such that $a_{k},\ b_{k}$ can be sufficiently large, then condition $2$ of Lemma \ref{PBfixed} can be easily verified. By Lemma \ref{starshape}, condition $1$ of Lemma \ref{PBfixed} also holds. In addition, we need to clarify $\mathrm{P_{1}}$ is an area-preserving mapping.
Indeed, the Poincar\'{e} map $\mathrm{P_{1}}: (z,w)\rightarrow\left(\overline{z}(2\pi,z,w),\overline{w}(2\pi,z,w)\right)$ can be expressed by
$$
\mathrm{P_{1}}=P^{1}\circ I^{1}\circ P^{0},
$$
where
$$
P^{0}:(z,w)\mapsto\left(\overline{z}(t_{1}-,z,w),\overline{w}(t_{1}-,z,w)\right),
$$
$$
I^{1}:\left(\overline{z}(t_{1}-),\overline{w}(t_{1}-)\right)\mapsto\left(\overline{z}(t_{1}-)+\Delta z(t_{1}),\overline{w}(t_{1}-)+\Delta w(t_{1})\right),
$$
$$
P^{1}:\left(\overline{z}(t_{1}+),\overline{w}(t_{1}+)\right)\mapsto\left(\overline{z}(2\pi),\overline{z}(2\pi)\right),
$$
and $\Delta z(t_{1}), \Delta w(t_{1})$ are given in \eqref{auxiliarytwoeq}.
Since $P^{j},\ j=0,1$ are symplectic by equation $z''+g(z)=0$ being conservative and $I^{1}$ is area-preserving with $\det(I^{1})=1$, the Poincar\'{e} map $\mathrm{P_{1}}$ is area-preserving.

Therefore, we can apply Lemma \ref{PBfixed} to ensure the existence of at least two fixed points of $\mathrm{P_{1}}$ in $\mathcal{A}_{k}$ ($k\geq n_{0}$). This means that \eqref{maineqautonomic} has at least two $2\pi$-periodic solutions with initial points in $\mathcal{A}_{k}$. Similarly, we can prove that $\mathrm{P_{1}}$ has at least two fixed points in $\mathcal{B}_{k}$ which correspond to two $2\pi$-periodic solutions of \eqref{maineqautonomic}. Since each periodic solutions of \eqref{maineqautonomic} is bounded by $\Gamma_{a_{k}}$ and $\Gamma_{b_{k}}$, then the above specified $2\pi$-periodic solutions of \eqref{maineqautonomic} constitute an infinite class.

The proof of Theorem \ref{mainresultautonomic} is thus completed.

\section{Nonautonomous Duffing impulsive equation}

We write \eqref{maineq} as an equivalent system of the form
\begin{equation}\label{maintwoeq}
\left\{\begin{array}{ll}
 x'=y,\quad y'=-g(x)+p(t),\quad t\neq t_{j};\\[0.2cm]
 \Delta x|_{t=t_{j}}=0,\\[0.2cm]
 \Delta y|_{t=t_{j}}=-2y\left(t_{j}-\right),\quad\ \ \ \ \ j=\pm1,\ \pm2,\cdots.\\
 \end{array}
 \right.
 \end{equation}
Let $\left(\overline{x}(t,x,y),\overline{y}(t,x,y)\right)$ be the solution of \eqref{maintwoeq} with the initial point $\left(\overline{x}(0),\overline{y}(0)\right)=(x,y)$. Also it follows from the condition $(H_{1})$ that every such solution exists on the whole $t$-axis (see \cite{Bainov93}). Then the Poincar\'{e} map $\mathrm{P_{2}}: \mathbb{R}^{2}\rightarrow \mathbb{R}^{2}$ is well defined by
$$(x,y)\rightarrow\left(\overline{x}(2\pi,x,y),\overline{y}(2\pi,x,y)\right).$$

By applying the transform $\overline{x}(t)=\gamma(t)\cos\theta(t), \overline{y}(t)=\gamma(t)\sin\theta(t)$ to \eqref{maintwoeq}, we get the equations for $\gamma(t)$ and $\theta(t)$,
\begin{equation}\label{maintwopolareq}
\left\{\begin{array}{ll}
 \gamma'=\gamma\cos\theta\cdot\sin\theta-g(\gamma\cos\theta)\sin\theta+p(t)\sin\theta,\\[0.2cm]
 \theta'=-\sin^{2}\theta-\displaystyle\frac{1}{\gamma}\left[g(\gamma\cos\theta)\cos\theta-p(t)\cos\theta\right],\quad t\neq t_{j};\\[0.4cm]
 \gamma(t_{j}+)=\gamma(t_{j}-),\\[0.2cm]
 \theta(t_{j}+)=\arctan(-\tan\theta (t_{j}-)),\ j=\pm1,\ \pm2,\cdots.\\
 \end{array}
 \right.
 \end{equation}

Let $\left(\overline{\gamma}(t,\gamma,\theta),\overline{\theta}(t,\gamma,\theta)\right)$ be the solution of \eqref{maintwopolareq} through the initial point $\left(\overline{\gamma}(0),\overline{\theta}(0)\right)=(\gamma,\theta)$. Then the map $\mathrm{P_{2}}$ also can be written in the polar coordinate form
\begin{equation}\label{p mapping form}
\gamma^{\ast}=\overline{\gamma}(2\pi,\gamma,\theta),\quad \quad \theta^{\ast}=\overline{\theta}(2\pi,\gamma,\theta)+2l\pi,
\end{equation}
where $l$ is an arbitrary integer. It can be easily seen that if $(\gamma,\theta)$ satisfies
\begin{equation}\label{gamma positive}
\overline{\gamma}(t,\gamma,\theta)>0, \quad \quad t\in[0,2\pi],
\end{equation}
then $\overline{\theta}(2\pi,\gamma,\theta)$ is well defined and continuous in $(\gamma,\theta)$, moreover,
$$\overline{\theta}(2\pi,\gamma,\theta+2\pi)=\overline{\theta}(2\pi,\gamma,\theta)+2\pi.$$
%%%%%%%%%%%%%%%%%%%%%%%%%%%%%%%%%%%%%%%%%%%%%%%%%%%%%%%%%%%%%%%%%%%%%%%%%%%%%%%%%%%%%%%%%%%%%%%%%%%%%%%%%%%%

Based on Lemma \ref{twistlemma} in Section 3, we first will study the intimate relation between equations \eqref{maineq} and \eqref{maineqautonomic}, especially the motions of their arguments under the polar coordinates. Let
$$
\Theta(\gamma,\theta)=\overline{\theta}(2\pi,\gamma,\theta)-\theta.
$$
The following Lemma will estimate the difference between $\Theta(\gamma,\theta)$ and $\Phi(\gamma,\theta)$ with the same initial point.
\begin{lemma}\label{sufficiently small}
For any $\varepsilon>0$, there exists $\gamma^{\ast}>0$ such that, for $\gamma\geq\gamma^{\ast}$,
$$
|\Theta(\gamma,\theta)-\Phi(\gamma,\theta)|=|\overline{\theta}(2\pi,\gamma,\theta)-\overline{\varphi}(2\pi,\gamma,\theta)|<\varepsilon.
$$
\end{lemma}

\Proof For $0<\varepsilon<\pi$, let $\left(\overline{z}(t,x,y),\overline{w}(t,x,y)\right)$ be the solution of \eqref{auxiliarytwoeq} with the initial point $\left(\overline{z}(0),\overline{w}(0)\right)=(x,y)$. Denote
$$
\left\{
\begin{array}{ll}
u(t)=u(t,x,y)=\overline{x}(t,x,y)-\overline{z}(t,x,y),\\[0.2cm]
v(t)=v(t,x,y)=\overline{y}(t,x,y)-\overline{w}(t,x,y).\\
\end{array}
 \right.
$$
Then we have
$$
\frac{du}{dt}=v,\quad\quad \frac{dv}{dt}=p(t)-g'(\vartheta(t))u(t),\quad t\neq t_{j},\ j=\pm1,\ \pm2,\cdots,
$$
where $\vartheta(t)=\overline{x}(t)+\lambda(t)\left(\overline{z}(t)-\overline{x}(t)\right)$, $0\leq\lambda(t)\leq1$.
Since the impulsive functions in \eqref{auxiliarytwoeq} and \eqref{maintwoeq} are same and linear, then
\begin{equation}\label{impulsive condition}
u(t_{j}+)=u(t_{j}),\quad\quad v(t_{j}+)=-v(t_{j}),\ j=\pm1,\ \pm2,\cdots.
\end{equation}
Let $\eta(t)=\left(u^{2}(t)+v^{2}(t)\right)^{\frac{1}{2}}$. Then for $t\neq t_{j}$, we have
$$
\eta\frac{d\eta}{dt}=uv+p(t)v-g'(\vartheta(t))uv.
$$
It follows from $(H_{1})$ that
\begin{equation}\label{derivative eta}
\Big|\frac{d\eta}{dt}\Big|\leq\frac{1}{2}(1+K)\eta+B,
\end{equation}
where $B=\displaystyle\sup_{t\in[0,2\pi]}|p(t)|$.

First we consider the difference between the arguments on $[0,t_{1}]$. The differential inequality \eqref{derivative eta} together with $\eta(0)=0$ yields
$$
\eta(t)\leq\frac{2B}{1+K}\left[e^{(1+K)\pi}-1\right]\equiv H_{0}
$$
for $t\in[0,t_{1}]$.
Denote
$$
\psi(t)=\psi(t,\gamma,\theta)=\overline{\theta}(t,\gamma,\theta)-\overline{\varphi}(t,\gamma,\theta),
$$
where $(\gamma,\theta)$ is the polar coordinate of $(x,y)$, that is, $(\gamma\cos\theta,\gamma\sin\theta)=(x,y)$. It is clear that if $\mid\psi(t)\mid<\pi$, then $\psi(t)$ is just the angle between the vectors $\left(\overline{x}(t),\overline{y}(t)\right)$ and $\left(\overline{z}(t),\overline{w}(t)\right)$. By the law of cosines, we have
$$
\cos\psi(t)=\displaystyle\frac{\overline{\gamma}^{2}(t)+\overline{\rho}^{2}(t)-\overline{\eta}^{2}(t)}{2\overline{\gamma}(t)\overline{\rho}(t)}\geq 1-\displaystyle\frac{H_{0}^{2}}{2\overline{\gamma}(t)\overline{\rho}(t)},\quad\quad t\in[0,t_{1}].
$$

On the other hand, $\overline{\gamma}(t)\geq\overline{\rho}(t)-\eta(t)\geq\overline{\rho}(t)-H_{0}$. Therefore, under the assumption that $|\psi(t)|<\pi$ and $\overline{\rho}(t)-H_{0}>0$, we have
\begin{equation}\label{cospsi}
\cos\psi(t)\geq1-\displaystyle\frac{H_{0}^{2}}{2\overline{\rho}(t)\left(\overline{\rho}(t)-H_{0}\right)},\quad\quad t\in[0,t_{1}].
\end{equation}
Note that
$$
F(\gamma,\theta)=\displaystyle\inf_{t\in[0,2\pi]}\overline{\rho}(t)=\displaystyle\inf_{t\in[0,2\pi]}\overline{\rho}(t,\gamma,\theta)
$$
becomes arbitrarily large if $\gamma$ is sufficiently large (see \cite{Qian15}). Then there is a constant $\gamma_{0}>0$ such that, for $\gamma\geq\gamma_{0}$ and $t\in[0,t_{1}]$,
\begin{equation}\label{cosdelta}
\overline{\rho}(t)-H_{0}>0,\quad\quad \displaystyle\frac{H_{0}^{2}}{2\overline{\rho}(t)\left(\overline{\rho}(t)-H_{0}\right)}<1-\cos\delta,
\end{equation}
where $\delta=\min\{\frac{\pi}{2},\frac{\varepsilon}{2}\}$. From \eqref{cospsi} and \eqref{cosdelta}, we conclude that if $|\psi(t)|<\pi$ on $[0,t_{1}]$, then the inequality
$$
\cos\psi(t)>\cos\delta
$$
holds, which implies
$$
|\psi(t)|<\delta,\quad\quad t\in[0,t_{1}].
$$

Now we just need to verify
\begin{equation}\label{psi1}
|\psi(t)|<\pi,\quad\quad t\in[0,t_{1}].
\end{equation}
Since $\psi(0)=0$ and $\psi(t)$ is continuous on $[0,t_{1}]$, there exists $a>0$ such that
\begin{equation}\label{psi2}
|\psi(t)|<\pi,\quad\quad t\in[0,a).
\end{equation}
To prove \eqref{psi1}, we need to verify $a> t_{1}$.
If it is not true, then the $a$ given above belongs to $(0,t_{1}]$ such that \eqref{psi2} holds and $|\psi(a)|=\pi$. Note that  inequality \eqref{cospsi}, together with \eqref{cosdelta} and \eqref{psi2} implies
$$
|\psi(t)|<\delta,\quad\quad t\in[0,a).
$$
Thus we obtain $|\psi(a)|\leq\delta\leq\frac{\pi}{2}<\pi$, which contradicts with $|\psi(a)|=\pi$, hence \eqref{psi1} holds and then for $t\in[0,t_{1}]$,
\begin{equation}\label{psismall}
|\psi(t)|<\delta\leq\frac{\varepsilon}{2}<\varepsilon.
\end{equation}

Next we turn to consider $\psi(t)$ on $(t_{1},2\pi]$. By \eqref{impulsive condition}, it follows that
$$
\eta(t_{1}+)=\eta(t_{1})\triangleq\eta_{1}.
$$ Integrating the inequality \eqref{derivative eta} on $(t_{1},2\pi]$, we have
$$
\eta(t)\leq\bigg|\frac{2B_{1}}{1+K}\left[e^{(1+K)\pi}-\frac{B}{B_{1}}\right]\bigg|\equiv H_{1},\quad\quad t\in(t_{1},2\pi],
$$
where $B_{1}=\frac{1+K}{2}\eta_{1}+B$. In the same way with the discussion on $[0,t_{1}]$, we obtain that there is a constant $\gamma_{1}>0$, such that for $\gamma\geq\gamma_{1}$ and $t\in (t_{1},2\pi]$,
\begin{equation}\label{cospsi1}
\cos\psi(t)\geq1-\displaystyle\frac{H_{1}^{2}}{2\overline{\rho}(t)\left(\overline{\rho}(t)-H_{1}\right)},
\end{equation}
and
\begin{equation}\label{cosdelta1}
\overline{\rho}(t)-H_{1}>0,\quad\quad \displaystyle\frac{H_{1}^{2}}{2\overline{\rho}(t)\left(\overline{\rho}(t)-H_{1}\right)}<1-\cos\delta,
\end{equation}
where $\delta$ is the same as before, provided $|\psi(t)|<\pi$ on $(t_{1},2\pi]$. It follows that if $|\psi(t)|<\pi$ on $(t_{1},2\pi]$, then the inequality
$$
|\psi(t)|<\delta
$$
holds.

Similarly, we just need to prove
\begin{equation}\label{psi21}
|\psi(t)|<\pi,\quad\quad t\in(t_{1},2\pi].
\end{equation}
By \eqref{psismall}, $\psi(t_{1})=\psi(t_{1}-)<\varepsilon$. Since the impulsive functions in \eqref{maintwopolareq} and \eqref{auxiliarytwopolareq} are same and $\psi(t)$ is continuous in $t_{1}$, $\psi(t_{1}+)$ can be small enough. Then there exists $a>0$ such that
\begin{equation}\label{psi22}
|\psi(t)|<\pi,\quad\quad t\in(t_{1},a).
\end{equation}
We just need to verify $a> 2\pi$.
Otherwise, there is an $a\in(t_{1},2\pi]$ such that \eqref{psi22} holds and $|\psi(a)|=\pi$. Note that  inequality \eqref{psi22}, together with \eqref{cospsi1} and \eqref{cosdelta1} implies
$$
|\psi(t)|<\delta,\quad\quad t\in(t_{1},a).
$$
Thus we obtain $|\psi(a)|\leq\delta\leq\frac{\pi}{2}<\pi$, which contradicts with $|\psi(a)|=\pi$. Then \eqref{psi21} holds and for $t\in(t_{1},2\pi]$,
$$
|\psi(t)|<\delta\leq\frac{\varepsilon}{2}<\varepsilon.
$$
In particular, $|\psi(2\pi)|<\varepsilon$.

Consequently, choosing $\gamma^{\ast}=\max\{\gamma_{0},\gamma_{1}\}$, for $\gamma\geq\gamma^{\ast}$, we have
$$
|\Theta(\gamma,\theta)-\Phi(\gamma,\theta)|=|\overline{\theta}(2\pi,\gamma,\theta)-\overline{\varphi}(2\pi,\gamma,\theta)|=|\psi(2\pi)|<\varepsilon.
$$
The proof of Lemma \ref{sufficiently small} is then completed.\qed

By Lemma \ref{twistlemma}, together with the small difference between $\overline{\theta}(2\pi)$ and $\overline{\varphi}(2\pi)$,
we can obtain the existence of $2\pi$-periodic solutions. However, according to Remark \ref{remark}, we cannot guarantee the application of the twist theorem for all sufficiently large $k$. Then the number of $2\pi$-periodic solutions of \eqref{maineq} is limited instead of infinity.\\

\noindent {\bf Proof of Theorem \ref{mainresult}} Let $c_{1}\geq c_{0}$ be sufficiently large, such that for $c\geq c_{1}$, $(\gamma,\theta)\in\Gamma_{c}$ implies $\gamma\geq\gamma^{\ast}$, where $\gamma^{\ast}$ is specified in Lemma \ref{sufficiently small}. There is no loss of generality to assume that $a_{k}, b_{k}\geq c_{1}$ for $k\geq n_{0}$, where $n_{0}$ is large enough. It follows that
\begin{equation}\label{gamma positive two}
\left\{
\begin{array}{ll}
\overline{\gamma}(t,\gamma,\theta)\geq\overline{\rho}(t)-H_{0}>0,\quad\quad t\in[0,t_{1}],\\[0.2cm]
\overline{\gamma}(t,\gamma,\theta)\geq\overline{\rho}(t)-H_{1}>0,\quad\quad t\in(t_{1},2\pi],\\
\end{array}
 \right.
\end{equation}
provided that $(\gamma,\theta)\in\mathcal{A}_{k}$ for $k\geq n_{0}$.
Thus the restriction $\mathrm{P_{2}}|\mathcal{A}_{k}$ can be written in \eqref{p mapping form}, where we put the integer $l=m$. Rewritten \eqref{p mapping form} in the following form
$$
\gamma^{\ast}=\overline{\gamma}(2\pi,\gamma,\theta),\quad \quad \theta^{\ast}=\theta+\Theta_{1}(\gamma,\theta),
$$
with $\Theta_{1}(\gamma,\theta)=\Theta(\gamma,\theta)+2m\pi$.

From Remark \ref{remark}, we first choose arbitrary finitely many integers $k_{1}$, $k_{2}$, $\cdots,$ $k_{n}$ with $k_{i}\geq n_{0},\ i=1,2,\cdots n$. Correspondingly, for each $k_{i}$, denote the curve by $\Gamma_{a_{k_{i}}},\Gamma_{b_{k_{i}}}$ and the annulus by $\mathcal{A}_{k_{i}}, \mathcal{B}_{k_{i}}$. To be precise, in Lemma \ref{twistlemma}, we denote $\beta_{2}=\beta_{2}^{i}$ for each $k_{i}$. Let $\varepsilon_{0}=\min\{\beta_{1}, \beta_{2}^{1},\beta_{2}^{2}\cdots,\beta_{2}^{n}\}$. Since $\beta_{1}>0$ and $\beta_{2}^{i}>0,\ i=1,2,\cdots n$, then $\varepsilon_{0}>0$. By Lemma \ref{sufficiently small}, choosing $\varepsilon=\varepsilon_{0}$, we obtain
$$
|\Theta_{1}(\gamma,\theta)-\Phi_{1}(\gamma,\theta)-2m\pi|<\varepsilon_{0},
$$
which together with \eqref{twistequaiton} yields
$$
\left\{
\begin{array}{ll}
\Theta_{1}(\gamma,\theta)<0,\quad\quad(\rho\cos\varphi,\rho\sin\varphi)\in \Gamma_{a_{k_{i}}};\\[0.2cm]
\Theta_{1}(\gamma,\theta)>0,\quad\quad(\rho\cos\varphi,\rho\sin\varphi)\in \Gamma_{b_{k_{i}}}.\\
\end{array}
 \right.
$$

This proves the validity of condition $3$ of Lemma \ref{PBfixed} for the restriction $\mathrm{P_{2}}|\mathcal{A}_{k_{i}}$ ($k_{i}\geq n_{0},\ i=1,2,\cdots n$). Since \eqref{gamma positive two} holds, then condition $2$ of Lemma \ref{PBfixed} can be easily verified. By Lemma \ref{starshape}, condition $1$ of Lemma \ref{PBfixed} also holds. Similar to the proof of Theorem \ref{mainresultautonomic}, area-preserving property of $\mathrm{P_{2}}$ is obtained.

Therefore, we can apply Lemma \ref{PBfixed} to ensure the existence of at least two fixed points of $\mathrm{P_{2}}$ in $\mathcal{A}_{k_{i}}$ ($k_{i}\geq n_{0},\ i=1,2,\cdots n$). This means that \eqref{maineq} has at least two $2\pi$-periodic solutions with initial points in $\mathcal{A}_{k_{i}}$ ($k_{i}\geq n_{0},\ i=1,2,\cdots n$). Similarly, we can prove that $\mathrm{P_{2}}$ has at least two fixed points in $\mathcal{B}_{k_{i}}$ ($k_{i}\geq n_{0},\ i=1,2,\cdots n$) which correspond to two $2\pi$-periodic solutions of \eqref{maineq}. Since each period solutions of \eqref{maineq} is bounded by $\Gamma_{a_{k_{i}}}, \Gamma_{b_{k_{i}}}$, then the above specified $2\pi$-periodic solutions of \eqref{maineq} constitute a finite class of solutions.

The proof of Theorem \ref{mainresult} is thus completed.

\section*{References}
\bibliographystyle{elsarticle-num}

\end{document}